\documentclass[twocolumn,letterpaper]{IEEEAerospaceCLS}

\usepackage[]{graphicx}
\usepackage[colorlinks=true, urlcolor=blue, linkcolor=black, citecolor=black]{hyperref}
\usepackage{amsmath}
\usepackage{amssymb}
\usepackage{float}
\usepackage{array}
\usepackage[version=4]{mhchem}
\usepackage{siunitx}
\usepackage{subcaption}
\usepackage{longtable,tabularx}
\usepackage{booktabs}
\usepackage[edges]{forest}
\usepackage{quoting}
\newcommand{\ignore}[1]{}
\begin{document}

\title{PLUTO: An Agentic AI Tool for Interactive Spacecraft Rendezvous Trajectory Design}

\author{%
Eleanor Brosius\textsuperscript{1},
Yuji Takubo\textsuperscript{1},
Daniele Gammelli\textsuperscript{1,2},
Simone D'Amico\textsuperscript{1},
Marco Pavone\textsuperscript{1,3}\\[6pt]
\textsuperscript{1}Department of Aeronautics and Astronautics, Stanford University\\
496 Lomita Mall, Stanford, CA 94305\\[4pt]
\textsuperscript{2}Italian Institute of Artificial Intelligence (AI4I), Turin, Italy 10128\\[4pt]
\textsuperscript{3}NVIDIA, Santa Clara, CA 95051\\[6pt]
\{elbro, ytakubo, gammelli, damicos, pavone\}@stanford.edu
\thanks{\footnotesize 979-8-3195-1048-8/27/$\$31.00$ \copyright2027 IEEE}
}

\maketitle
\thispagestyle{plain}
\pagestyle{plain}

\begin{abstract}
Modern space missions require trajectory design methods capable of efficiently adapting to diverse mission objectives and operational constraints. Despite the maturity of non-convex trajectory optimization techniques, their broader adoption remains limited by the substantial domain expertise required to tailor problem formulations across varying mission scenarios. Recent advances in large language model (LLM) reasoning and agentic AI coding capabilities offer a promising pathway to improving the accessibility of advanced trajectory optimization techniques. This paper introduces \textbf{\textsc{PLUTO}} (\textbf{P}lain \textbf{L}anguage \textbf{U}nderstanding for \textbf{T}rajectory \textbf{O}ptimization), an interactive trajectory design framework that integrates agentic AI coding tools within a sequential convex programming (SCP) architecture. The proposed system maps high-level semantic inputs to well-defined mathematical constraints while preserving a structured representation suitable for SCP. These constraints are processed through an auto-convexification pipeline and incorporated directly into executable optimal control formulations, thereby systematically connecting user intent and trajectory generation. Visual constraint verification and scenario parameterization further support an intuitive environment for trajectory design across a wide range of application scenarios. Within an evaluation of 50 natural-language rendezvous mission prompts spanning 10 constraint types, \textsc{PLUTO} generated a trajectory for every prompt, and 94\% of the resulting trajectories satisfied both the quantitative and qualitative requirements of the mission intent. Ultimately, this work explores the potential of agentic AI to bridge the gap between advanced optimization algorithms and engineer-focused trajectory design workflows.

\vspace{-5.0mm}
\end{abstract}

\section{Introduction}

\begin{figure}[t]
    \centering
    \includegraphics[width=1\linewidth]{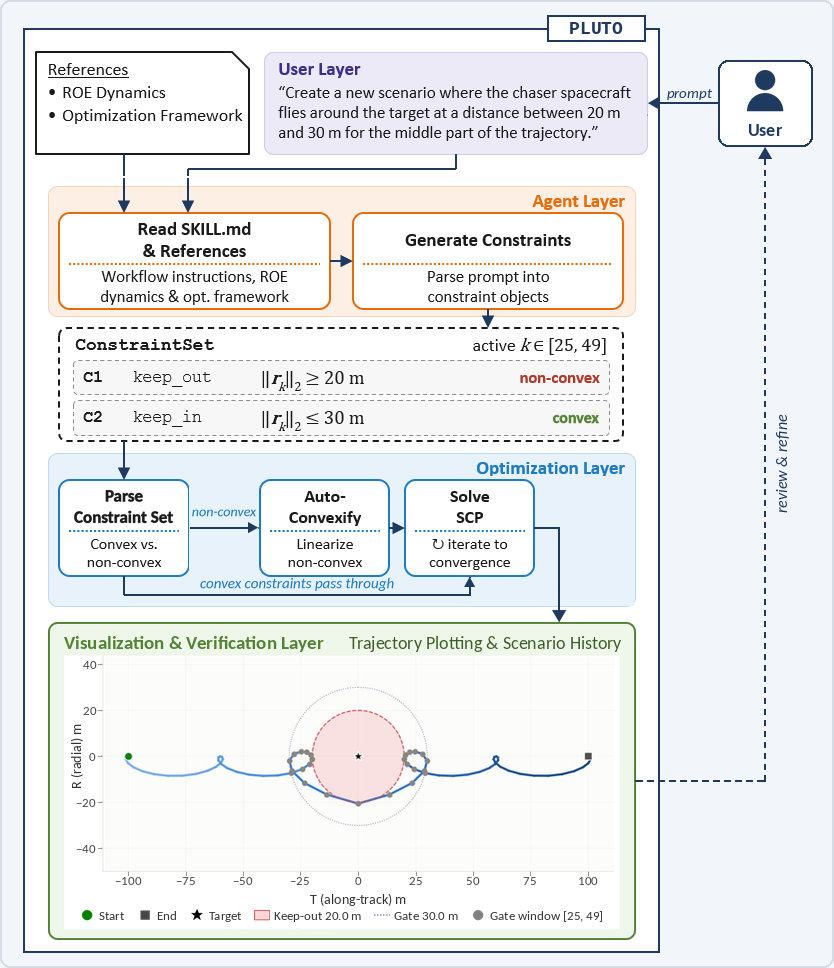}
    \caption{PLUTO (Plain Language Understanding for Trajectory Optimization) architecture. The agent parses a natural-language prompt into
convexity-classified constraint objects, which the optimization layer sorts,
convexifies, and solves via SCP. Results are displayed in the dashboard for user
review and iteration.}
    \label{fig:hero}
\end{figure}

Agentic artificial intelligence (AI) is creating new opportunities for automating complex scientific and engineering workflows. 
Built on large multimodal models trained on internet-scale data, these systems combine reasoning and code generation with tool use, memory, and iterative feedback to carry out multistep tasks. 
Emerging applications across a range of disciplines demonstrate their ability to coordinate specialized software, incorporate domain knowledge, and use computational results to guide subsequent actions across a broad range of tasks \cite{wang2024chatstl,ahmaditeshnizi2024optimus}. 
These capabilities are particularly promising for engineering design, where translating high-level objectives into computational models often requires substantial expertise and repeated iteration.

Spacecraft trajectory design places heavy demands on this expertise. 
Mission-level objectives and operational requirements must be translated into mathematical formulations that can be solved and validated computationally. 
These formulations must capture the relevant dynamics and constraints while remaining adaptable as mission concepts evolve.
Engineers must therefore repeatedly formulate, solve, and validate optimal control problems (OCPs) as mission concepts are refined.

General-purpose methods for non-convex optimization provide a flexible foundation for a broad class of trajectory design problems. 
In particular, approaches such as SCP offer an effective framework for handling the non-convexities that arise in many practical formulations \cite{malyuta_scp_2022}. 
However, adapting these methods to new or evolving mission scenarios still requires substantial expert effort to construct mathematically consistent formulations and translate them into executable implementations.

Automating this process requires more than general language understanding or code generation. 
Mission requirements must be translated into consistent dynamical models and optimization formulations that respect solver assumptions and numerical considerations, and the resulting solutions must then be assessed for both numerical feasibility and fidelity to the intended mission. 
Because these decisions span astrodynamics, optimal control, numerical optimization, and solution verification, reliable agentic automation requires a domain-specific scaffold that provides explicit modeling guidance, interfaces to established optimization tools, and procedures for inspecting and validating intermediate outputs.

This paper introduces \textbf{\textsc{PLUTO}}
(\textbf{P}lain \textbf{L}anguage\break
\textbf{U}nderstanding for \textbf{T}rajectory \textbf{O}ptimization), an interactive framework for agentic AI-assisted spacecraft trajectory design. 
As illustrated in Fig.~\ref{fig:hero}, \textsc{PLUTO} guides an agent through the translation of natural-language mission requirements into structured constraint definitions and executable OCP formulations. 
A dedicated SCP backend incorporates convex constraints directly and processes supported non-convex constraints through an auto-convexification pipeline. 
Compiled domain knowledge and workflow guidance support consistent modeling decisions, while plotting and dashboard utilities expose trajectory behavior and constraint satisfaction for user review and iteration. 
When processed correctly, this separation between agent-driven formulation and numerical optimization ensures that the resulting trajectory satisfies the prescribed constraints by construction, a property that is particularly critical in aerospace applications.

The contributions of this work are as follows:
\begin{itemize}	
    \item The introduction of \textsc{PLUTO}, an agentic AI framework for spacecraft trajectory optimization that integrates requirement interpretation, OCP formulation, numerical solution, and verification into a unified interactive workflow.
    
	\item A suite of domain-specific tools supporting \textsc{PLUTO}, including SCP and auto-convexification, agent workflows and guardrails, compiled dynamics and constraint knowledge, and utilities for trajectory visualization, constraint inspection, and dashboard generation.
    
	\item An empirical evaluation across 50 spacecraft rendezvous scenarios, in which \textsc{PLUTO} generated a trajectory for 100\% of prompts and 94\% of the resulting trajectories satisfied both numerical constraint checks and semantic consistency with the original mission intent, together with an analysis of failure modes across the multilayer framework.
\end{itemize}

\section{Related Work}

The proposed framework builds upon prior advances in numerical trajectory optimization, learning-based trajectory generation, and LLM-based approaches for synthesizing mathematical constraints from natural-language specifications. This section provides a brief overview of these research areas and summarizes the key limitations that motivate the proposed approach.

\subsection{Spacecraft Trajectory Optimization and Learning-Based Trajectory Generation}

Trajectory optimization serves as a foundational tool for both onboard guidance and the offline design of complex space missions, as it enables the simultaneous treatment of formal constraint satisfaction and control-cost optimality. To address the non-convexity arising from nonlinear dynamics and operational constraints, a variety of optimization methods have been developed, including SCP \cite{malyuta_scp_2022}, together with theoretical analyses establishing convergence properties and guarantees \cite{mao2016successive}, \cite{oguri2023successive}, \cite{bonalli_2019_gusto}. Because the performance of non-convex trajectory optimization can remain highly sensitive to the initial guess, intelligent warm-starting strategies have also been investigated, including learning-based approaches for generating high-quality initial trajectories \cite{Banerjee_2020}, \cite{art_ieeeaero24}, \cite{briden2025transformer}, \cite{kim2022guided}.

More recently, with the emergence of large language models (LLMs), their integration into spacecraft guidance, navigation, and control architectures has received increasing attention. Some studies investigate LLMs for direct reasoning and decision-making \cite{foutter2025space-llava}, \cite{rodriguez2024language}, \cite{jain2026autonomous}, while others consider language-conditioned trajectory generation for warm-starting SCP \cite{takubo2026language} or scaffolded LLM frameworks that systematically construct mission configurations within a predefined design space before invoking trajectory optimization \cite{takubo2026orbittamp}. 
A related tooling-based agentic framework for spacecraft formation control is presented in \cite{zhang2025llm}. 
However, these approaches largely operate over fixed optimization formulations, with task-specific constraints either predefined by the system designer or not explicitly incorporated into the optimization problem.

This leaves a key gap between natural-language mission intent and formal trajectory-optimization models. 
Our work addresses this gap by enabling zero-shot synthesis of mathematical constraints from natural-language specifications and incorporating those constraints directly into the optimization problem. 
In this way, operator intent can modify the formulation at the constraint level, rather than only selecting among predefined configurations or providing inputs to an otherwise fixed optimizer.

\subsection{Large Language Models for Constrained Optimization}

Recent work in optimization and operations research has investigated the translation of natural-language descriptions into structured mathematical formulations \cite{ramamonjison2023nl4opt}. 
Building on this direction, agentic systems have been developed to construct optimization problems, generate executable solver code, and iteratively refine combinatorial optimization formulations \cite{ahmaditeshnizi2024optimus}. Related efforts have addressed the translation of language into temporal-logic specifications  \cite{wang2024chatstl,chen2023nl2tl,van2024vernacopter}, the enforcement of language-specified constraints through robotic safety mechanisms \cite{yang2024plug}, \cite{feng2025words}, and the automated reformulation of non-convex optimization problems \cite{peng2026nc2c}. 
More closely related to optimal control, NARRATE combines sequential task planning with the generation of objective and constraint functions for model predictive control \cite{ismail2024narrate}.

Building on these developments and advances in language-driven code generation \cite{liang2023code}, this work brings constraint synthesis directly into an SCP-based spacecraft trajectory optimization workflow. 
The agent generates executable, automatically differentiable constraint functions that can be incorporated into the optimization problem and, when non-convex, linearized within successive convex subproblems. This closes the loop between natural-language mission intent and numerical spacecraft trajectory optimization, enabling task-specific operational constraints to be directly injected into the optimization formulation itself rather than being limited by a predefined constraint set.

\section{Problem Statement}

The central question considered in this work is whether an agentic AI system can translate high-level spacecraft mission intent into a mathematically valid, executable, and verifiable trajectory optimization problem. 
Rather than requiring the user to specify equations, constraint classes, or solver code, we seek an agent that can reason over a natural-language request, the current mission context, and domain-specific guidance to determine how the underlying optimal control problem should be modified.

Let $p$ denote a natural-language mission specification and let $\mathcal{S}$ denote the available scenario context, including the dynamics, state representation, reference frames, solver interfaces, and domain knowledge exposed to the agent. We seek an agentic mapping:
\begin{equation}
    \mathcal{A} :
    (p,\mathcal{S})
    \longmapsto
    \mathcal{C}
    =
    \{C_1,\ldots,C_M\},
    \label{eq:agent_mapping}
\end{equation}
where each $C_i$ is an executable constraint object defining its mathematical form, activation interval, convexity class, and numerical implementation.
The generated constraint set modifies a discrete-time optimal control problem of the form:

\begin{subequations}
\label{eq:ocp}
\begin{alignat}{2}
    \min_{\{{x}_k, {u}_k\}_{k=1}^{N}}
    \quad &
    \sum_{k=1}^{N-1}
    \mathcal{J}_k(\boldsymbol{x}_k,\boldsymbol{u}_k)
    + \mathcal{J}_N(\boldsymbol{x}_N)
    \label{eq:ocp_obj}
    \\[1mm]
    \text{subject to} \quad
    &
    \boldsymbol{x}_{k+1}
    =
    f_k(\boldsymbol{x}_k,\boldsymbol{u}_k),
    \quad
    && k \in \mathbb{Z}_{1:N-1},
    \label{eq:ocp_dyn}
    \\
    &
    \boldsymbol{x}_1 = \boldsymbol{x}_0,
    \qquad
    \boldsymbol{x}_N = \boldsymbol{x}_f,
    &&
    \label{eq:ocp_bc}
    \\
    &
    g_j(\boldsymbol{x}_k,\boldsymbol{u}_k) \leq 0,
    \quad
    && k \in \mathcal{K}_j,
    \quad j \in \mathcal{I},
    \label{eq:ocp_ineq}
    \\
    &
    h_j(\boldsymbol{x}_k,\boldsymbol{u}_k) = 0,
    \quad
    && k \in \mathcal{K}_j,
    \quad j \in \mathcal{E}.
    \label{eq:ocp_eq}
\end{alignat}
\end{subequations}
Here, $\boldsymbol{x}_k$ and $\boldsymbol{u}_k$ denote the system state and control at discrete step $k$, respectively, and $N$ the number of timesteps in the trajectory. 
The function $f_k(\boldsymbol{x}_k,\boldsymbol{u}_k)$ defines the discrete-time dynamics, while $\boldsymbol{x}_0$ and $\boldsymbol{x}_f$ specify the prescribed initial and terminal states. 
For the spacecraft rendezvous problems considered in this work, the control input is the impulsive maneuver, such that $\boldsymbol{u}_k = \Delta\boldsymbol{v}_k$. 
The functions $\mathcal{J}_k$ and $\mathcal{J}_N$ denote the stage and terminal costs, respectively. 
A common objective for spacecraft rendezvous is to minimize total maneuver effort as a proxy for fuel consumption, with $\mathcal{J}_k = \lVert \Delta\boldsymbol{v}_k \rVert_2$ and $\mathcal{J}_N = 0$, yielding the objective $\sum_{k=1}^{N-1}\lVert \Delta\boldsymbol{v}_k \rVert_2$. 
The set $\mathbb{Z}_{1:N-1}$ denotes the integer indices $\{1,\ldots,N-1\}$.

The functions $g_j$ and $h_j$ denote inequality and equality constraints, respectively, with $j \in \mathcal{I}$ for inequality constraints and $j \in \mathcal{E}$ for equality constraints.
For each constraint $j$, $\mathcal{K}_j$ denotes the set of trajectory steps over which it is enforced.
Thus, a constraint may be imposed over an interval $\mathcal{K}_j$, at a single step $k^\star$, or over the full trajectory, depending on the mission requirement. 
In the proposed framework, the agent may synthesize the functions $g_j$ and $h_j$, determine their activation sets $\mathcal{K}_j$, and generate the corresponding numerical implementations from the natural-language specification $p$.

Importantly, successful translation requires more than producing a plausible equation. 
The generated formulation must preserve the semantic intent of the request, respect the dynamics and coordinate conventions of the scenario, satisfy the mathematical requirements of the optimization method, and admit an executable implementation compatible with the numerical backend.
The resulting trajectory must then be checked against the original request, closing the loop between language interpretation, mathematical formulation, numerical optimization, and semantic verification.

\section{Methodology}

\textsc{PLUTO} operationalizes the language-to-optimization mapping introduced in the previous section through four layers: a user interface, an agent scaffold, a numerical optimization backend, and a visualization and verification layer.
The layers communicate through structured interfaces so that the agent is responsible for interpreting mission intent and constructing problem-specific artifacts, while the underlying numerical optimization routines remain explicitly defined in code.
This separation allows the agent to modify the trajectory optimization formulation without requiring it to implement or alter the solver itself.

At a high level, given a natural-language specification $p$ and scenario context $\mathcal{S}$, the agent produces a structured set of mathematical constraints and associated verification artifacts.
The constraints are passed to the optimization backend, incorporated into the OCP, and solved using SCP.
The resulting trajectory and constraint evaluations are then exposed to the user for inspection and subsequent iteration.
The following sections describe each layer and the interfaces between them.

\subsection{User Layer}
The user interacts with \textsc{PLUTO} through an agentic coding environment using a natural-language description of the desired spacecraft behavior.
The request may define a new trajectory design scenario or modify an existing problem without requiring the user to provide the corresponding mathematical constraints, convexity classifications, or solver implementations.

A \textsc{PLUTO} workspace contains the executable trajectory optimization code together with the contextual information required by the agent.
The initial user request directs the agent to a workflow controller, which serves as the entry point to the framework and specifies how the remaining resources and tools should be used.
In this way, the user-facing interface remains language-based while the downstream artifacts conform to the structured interfaces expected by the optimization backend.

\subsection{Agent Layer}
The agent layer translates the user request into an executable representation of the trajectory optimization problem.
Rather than allowing the agent to operate directly on the numerical solver without structure, \textsc{PLUTO} provides a domain-specific scaffold that constrains the workflow while retaining flexibility in how individual mission requirements are interpreted and formulated.

\textbf{Workflow Controller.} The workflow controller defines the high-level sequence of operations the agent must perform and the artifacts that must be produced before execution is complete.
It specifies guardrails on agent behavior, identifies the available modes of operation, and directs the agent to the domain-specific resources relevant to each stage of the workflow.
The controller therefore provides a common execution structure across agentic coding models while leaving the underlying model responsible for reasoning over the particular user request.

\textbf{Domain Knowledge and Reference Material.} The workflow controller provides access to three complementary forms of reference material.
\textit{Templates} define the structure of files and objects generated by the agent, ensuring compatibility with downstream interfaces.
\textit{Scenario documents} describe the relevant dynamics, geometry, coordinate conventions, and state representations.
Finally, \textit{formulation guides} provide instructions for decisions that require domain expertise, including the construction and numerical implementation of optimization constraints.

Together, these resources constitute the scenario context $\mathcal{S}$ introduced in the problem statement.
They allow the agent to ground natural-language requirements in the conventions of the trajectory design problem rather than relying only on the pretrained model's internal knowledge.

\textbf{Constraint Synthesis.} For each mission requirement identified in the user request, the agent determines an appropriate mathematical representation and the portion of the trajectory over which it should be enforced.
The resulting constraint is encoded directly as a structured Python object containing its mathematical implementation, activation set, and convexity classification.

The generated constraint objects are collected in a \texttt{ConstraintSet}, which provides the primary interface between the agent and optimization layers.
The agent additionally generates plotting and verification configurations associated with the synthesized requirements.
These artifacts preserve the connection between the original semantic request, its mathematical formulation, and the resulting trajectory.

\subsection{Optimization Layer}
\textsc{PLUTO} uses SCP as its numerical optimization backbone. 
This choice provides a common interface for incorporating a broad range of trajectory constraints while retaining the computational structure of convex optimization at each iteration. 
CVXPY is used to construct the convex subproblems, while JAX provides automatic differentiation for supported non-convex constraint functions.

\begin{figure}[t]
    \centering
    \includegraphics[width=1\linewidth]{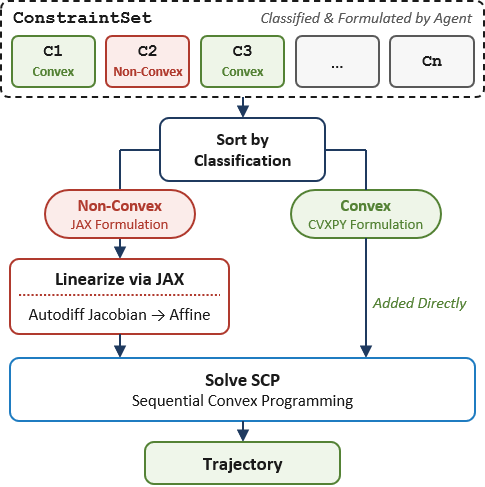}
    \caption{Constraint processing within the \textsc{PLUTO} optimization layer. Agent-generated constraints are routed according to convexity before being incorporated into the SCP subproblem.}
    \label{fig:convexifyflow}
\end{figure}

The optimization layer receives the agent-generated \penalty0\hspace{0pt}\mbox{\texttt{ConstraintSet}} and separates its members according to convexity.
Constraints represented directly in CVXPY are incorporated into the corresponding convex subproblem without further approximation.
Supported non-convex constraints represented in JAX are instead passed through the auto-convexification pipeline shown in Fig.~\ref{fig:convexifyflow}.

For a differentiable non-convex constraint:
\begin{equation}
    c(\boldsymbol{z}) \leq 0,
\end{equation}
where $\boldsymbol{z}$ collects the relevant state and control variables, a first-order approximation about the current SCP reference point $\bar{\boldsymbol{z}}$ is given by:
\begin{equation}
    c(\bar{\boldsymbol{z}})
    +
    \nabla c(\bar{\boldsymbol{z}})^\top
    \left(
    \boldsymbol{z}-\bar{\boldsymbol{z}}
    \right)
    \leq 0.
    \label{eq:constraint_linearization}
\end{equation}
The required derivatives are obtained automatically using JAX.
The resulting affine representation is converted into a CVXPY-compatible form and incorporated into the current convex subproblem.
As the SCP iterations progress, supported non-convex constraints are re-linearized about the updated reference trajectory.

\begin{table*}[t]
\centering
\caption{Constraint types used in the evaluation to form test prompts. $\boldsymbol{r}_k$ and $\boldsymbol{v}_k$ are the chaser's RTN position and velocity relative to the target at step $k$, $[\boldsymbol{r}_k;\boldsymbol{v}_k]=\Psi_k\boldsymbol{x}_k$ with $\boldsymbol{x}_k$ the relative orbital elements; $\Delta\boldsymbol{v}_k$ is the impulsive maneuver at step $k$; $\mathcal{K}$ is a constraint's active window and $k^\star$ a single step; $\hat{\boldsymbol{e}}_R,\hat{\boldsymbol{e}}_T,\hat{\boldsymbol{e}}_N$ are the RTN unit vectors. $\boldsymbol{c}$ is a zone center, $Q$ a rotation and $D$ the diagonal of semi-axes; $\hat{\boldsymbol{u}}$ is a cone or approach axis.}
\label{tab:constraint_types}
\renewcommand{\arraystretch}{1.55}
\begin{tabular}{lllll}
\hline
ID & Type & Constraint & Formulation & Convexity \\
\hline
T01 & Keep-out zone & Sphere or ellipsoid &
$\lVert D^{-1} Q^\top (\boldsymbol{r}_k - \boldsymbol{c}) \rVert_2 \ge 1, \; k \in \mathcal{K}$ \quad (sphere: $D = \rho I$) &
Non-convex \\
\hline
T02 & Keep-in region & Ball around the target &
$\lVert \boldsymbol{r}_k \rVert_2 \le \rho, \; k \in \mathcal{K}$ or $k = k^\star$ &
Convex \\
\hline
T03 & Waypoint & Pass a point &
$\lVert \boldsymbol{r}_{k^\star} - \boldsymbol{w} \rVert_2 \le \epsilon$ &
Convex \\
 & & Hold near a point &
$\lVert \boldsymbol{r}_k - \boldsymbol{w} \rVert_2 \le \epsilon, \; k \in \mathcal{K}$ &
Convex \\
\hline
T04 & Planar bounds & Half-space &
$\boldsymbol{n}^\top \boldsymbol{r}_k \ge d, \; k \in \mathcal{K}$ &
Convex \\
 & & Corridor &
$\boldsymbol{l} \le E\,\boldsymbol{r}_k \le \boldsymbol{u}, \; k \in \mathcal{K}$ &
Convex \\
\hline
T05 & Approach cone & Cone &
$\lVert (I - \hat{\boldsymbol{u}}\hat{\boldsymbol{u}}^\top) \boldsymbol{r}_k \rVert_2 \le \tan\theta \; \hat{\boldsymbol{u}}^\top \boldsymbol{r}_k, \; k \in \mathcal{K}$ &
Convex\\
\hline
T06 & Speed limit & Speed &
$\lVert \boldsymbol{v}_k \rVert_2 \le v_{\max}, \; k \in \mathcal{K}$ &
Convex \\
 & & Velocity component &
$\lvert \hat{\boldsymbol{e}}_i^\top \boldsymbol{v}_k \rvert \le v_{\max}, \; k \in \mathcal{K}$ &
Convex  \\
\hline
T07 & Maneuver magnitude & Cap &
$\lVert \Delta\boldsymbol{v}_k \rVert_2 \le u_{\max}, \; k \in \mathcal{K}$ &
Convex \\
 & & Minimum &
$\lVert \Delta\boldsymbol{v}_{k^\star} \rVert_2 \ge u_{\min}$ &
Non-convex \\
\hline
T08 & $\Delta v$ budget & Windowed budget &
$\sum_{k \in \mathcal{K}_j} \lVert \Delta\boldsymbol{v}_k \rVert_2 \le B$ \; for each window $\mathcal{K}_j$ &
Convex  \\
\hline
T09 & Thrust schedule & Coast window &
$\Delta\boldsymbol{v}_k = \boldsymbol{0}, \; k \in \mathcal{K}$ &
Convex  \\
 & & Maneuver count &
$\Delta\boldsymbol{v}_k = \boldsymbol{0}, \; k \notin \mathcal{K}_{\mathrm{burn}}, \;\; \lvert \mathcal{K}_{\mathrm{burn}} \rvert \le m$ chosen &
Convex \\
\hline
T10 & Misc & Range-dependent speed &
$\lVert \boldsymbol{v}_k \rVert_2 \le v_0 + \kappa \lVert \boldsymbol{r}_k \rVert_2, \; k \in \mathcal{K}$ &
Non-convex \\
 & & Glide slope &
$-\boldsymbol{r}_k^\top \boldsymbol{v}_k / \lVert \boldsymbol{r}_k \rVert_2 \le c_0 + \kappa \lVert \boldsymbol{r}_k \rVert_2, \; k \in \mathcal{K}$ &
Non-convex \\
 & & Out-of-plane exclusion &
$\sigma\,\hat{\boldsymbol{e}}_N^\top \boldsymbol{r}_k \ge d, \; \sigma \in \{+1, -1\}$ chosen, $k \in \mathcal{K}$ &
Convex  \\
\hline
\end{tabular}
\end{table*}

\subsection{Visualization and Verification Layer}
The final layer exposes the optimized trajectory and associated constraint information for analysis and user review.
For each synthesized requirement, the agent generates a plotting or verification configuration that is processed by a common visualization pipeline.
This allows scenario-specific quantities and constraint geometries to be displayed without requiring the agent to repeatedly implement low-level plotting functionality.

The resulting dashboard presents the spacecraft trajectory together with the active constraint regions and relevant optimization diagnostics.
Numerical verification evaluates whether the generated trajectory satisfies the implemented constraints over their prescribed activation sets.
Complementary semantic verification allows the user to inspect whether the resulting trajectory and displayed constraint geometry are consistent with the original natural-language request.
This distinction is important because a numerically satisfied constraint may nevertheless reflect an incorrect interpretation of the user's intent.

The dashboard is generated as an HTML artifact that records the outputs associated with each scenario.
By retaining the original request, generated formulation, optimized trajectory, and verification outputs within a common workflow, \textsc{PLUTO} supports iterative trajectory design in which the user can inspect the current solution and refine the mission specification in subsequent interactions.

\section{Results}

This section describes the experimental setup and evaluation metrics used to assess \textsc{PLUTO}, followed by the corresponding experimental results.

\subsection{Experimental Setup}
To evaluate PLUTO, we consider ten constraint types spanning convex and non-convex requirements relevant to spacecraft rendezvous. 
Each test scenario augments a common baseline OCP with one or more prompt-specific constraints. The baseline problem is given by:
\begin{subequations}
\label{eq:ocp_rendezvous}
\begin{align}
    \min_{\{\boldsymbol{x}_k, \boldsymbol{u}_k\}_{k=1}^{N}}
    \quad &
    \sum_{k=1}^{N-1} \|\boldsymbol{u}_k\|_2
    \label{eq:base_obj}
    \\[1mm]
    \text{subject to} \quad
    &
    \boldsymbol{x}_{k+1} = \Phi_k\boldsymbol{x}_k+\Gamma_k\boldsymbol{u}_k,
    && k \in \mathbb{Z}_{1:N-1},
    \label{eq:base_dyn}
    \\
    &
    \boldsymbol{x}_1 = \boldsymbol{x}_0,
    \label{eq:base_bc}
    \\
    &
    \|(\Psi_k \boldsymbol{x}_k)_{1:3}\|_2 \geq r_{\mathrm{KOZ}},
    && k \in \mathbb{Z}_{1:\lceil 0.8(N-1)\rceil}.
    \label{eq:base_koz}
\end{align}
\end{subequations}
where $\boldsymbol{x}_k \in \mathbb{R}^6$ denotes the Relative Orbital Elements (ROE) state~\cite{damico_phd_2010}, and $\boldsymbol{u}_k\in\mathbb{R}^3$ denotes an impulsive velocity change expressed in the target-centered radial--tangential--normal (RTN) frame. 
$\Phi_k = \Phi(t_{k+1}, t_k) \in \mathbb{R}^{6\times6}$ is the ROE state transition matrix with the secular $J_2$ effect~\cite{koenig2017new}, and $\Gamma_k = \Gamma(t_k) \in \mathbb{R}^{6\times3}$ is the control input matrix~\cite{damico_phd_2010}.
Finally, $\boldsymbol{\Psi}_k = \Psi(t_k)$ denotes the linear mapping from the ROE to the Cartesian position and velocity in the RTN frame \cite{damico_phd_2010}.
A spherical keep-out zone of radius $r_{\mathrm{KOZ}}$ is imposed around the target over the first 80\% of the trajectory. 
This common constraint prevents newly generated requirements from trivially altering the intended temporal structure of the rendezvous. All experiments use the same initial condition and time horizon to provide a consistent basis for comparison, while the remaining problem structure may be modified through the natural-language prompt.

The evaluation set consists of 50 prompts constructed from the ten constraint types summarized in Table~\ref{tab:constraint_types}, which also gives the mathematical formulation of each. 
Thirty prompts contain a single constraint type, ten combine two constraint types, and ten combine three. This distribution evaluates performance across increasing levels of compositional complexity while also exposing constraint classes that are more difficult for the agent to formulate correctly. 
Type T10, denoted \emph{Miscellaneous}, intentionally contains less standard and more challenging constraints.

Before evaluating \textsc{PLUTO}, a reference formulation was manually constructed and solved for each prompt. 
Each reference problem was required to be feasible within the underlying optimization framework, to produce a trajectory, and to satisfy the same numerical and semantic verification criteria later applied to the agent-generated solution.
This procedure establishes that each prompt corresponds to a realizable scenario and helps distinguish failures in agent-generated formulations from infeasibility arising from the underlying dynamics or mission geometry.

The prompts combine requirement-style specifications with more natural descriptions of desired spacecraft behavior, requiring the agent to interpret both explicit numerical requirements and higher-level semantic instructions. All experiments use \textit{claude-sonnet-5} as the agentic model, selected to provide a balance between formulation capability and computational cost. The present evaluation therefore focuses on the behavior of the complete \textsc{PLUTO} framework with a fixed frontier model. 
Evaluating model dependence, including the use of open-weight models, is left to future work.

\subsection{Evaluation Metrics}
\textsc{PLUTO} is evaluated using four binary metrics that capture performance at different stages of the language-to-trajectory pipeline.
For each prompt $i$, let
$M_i^{\mathrm{form}}$,
$M_i^{\mathrm{solve}}$,
$M_i^{\mathrm{opt}}$, and
$M_i^{\mathrm{verify}}$
denote the indicators for Constraint Formulation, Solver Completion, Optimal Convergence, and Prompt Verification, respectively.

\begin{itemize}
    \item \textbf{Constraint Formulation.} A prompt passes if the agent's formulation matches the reference formulation constructed for that prompt, meaning every required constraint is represented in its expected mathematical form and activation window, and no unintended constraints are introduced. Because natural language can admit more than one valid interpretation, a mismatch does not necessarily indicate an error. This metric instead measures how consistently the agent's interpretation aligns with the intended formulation. We define:
\begin{equation}
M_i^{\mathrm{form}} = \mathbf{1}\left[\begin{array}{c}\text{all constraints match the reference} \\ \text{and are activated as expected}\end{array}\right].
\end{equation}

    \item \textbf{Solver Completion.}
    A prompt passes if the generated formulation is successfully processed by the optimization pipeline and the solver returns a trajectory.
    This criterion does not require the solver to terminate under its optimal stopping condition.
    We define:
    \begin{equation}
        M_i^{\mathrm{solve}}
        =
        \mathbf{1}
        \left[
        \text{a trajectory is returned}
        \right].
        \label{eq:metric_completion}
    \end{equation}

    \item \textbf{Optimal Convergence.}
    A prompt passes if the SCP solver terminates after satisfying its prescribed optimality and feasibility conditions.
    Termination due to stagnation, an iteration limit, a penalty limit, or another non-optimal stopping condition is recorded as a failure for this metric, even if a trajectory is returned.
    We define:
    \begin{equation}
        M_i^{\mathrm{opt}}
        =
        \mathbf{1}
        \left[
        \begin{array}{c}
        \text{SCP terminates under its }
        \\
        \text{optimal stopping condition}
        \end{array}
        \right].
        \label{eq:metric_optimality}
    \end{equation}
    Because non-optimal termination may result from either the generated formulation or the numerical optimization process, this metric is interpreted jointly with Constraint Formulation.

    \item \textbf{Prompt Verification.}
    A prompt passes only if the resulting trajectory satisfies every requirement expressed in the original natural-language specification.
    Let $\mathcal{R}_i$ denote the set of requirements extracted from prompt $i$, and let $V_{i,r}\in\{0,1\}$ indicate whether requirement $r\in\mathcal{R}_i$ is satisfied.
    Quantitative requirements are evaluated using automated numerical checks over the relevant trajectory steps, including waypoint tolerances, separation distances, velocity limits, and maneuver bounds.
    Qualitative requirements that are not completely described by a scalar threshold are assessed from the generated trajectory and constraint visualizations.
    The prompt-level verification metric is therefore
    \begin{equation}
        M_i^{\mathrm{verify}}
        =
        \prod_{r\in\mathcal{R}_i} V_{i,r}
        =
        \mathbf{1}
        \left[
        V_{i,r}=1,\;\forall r\in\mathcal{R}_i
        \right].
        \label{eq:metric_verification}
    \end{equation}
\end{itemize}

Prompt Verification is the primary measure of end-to-end success because it directly evaluates whether the final trajectory satisfies the user's stated mission intent.
The remaining metrics provide diagnostic information about intermediate stages of the pipeline and help distinguish failures in constraint synthesis from those arising during numerical optimization.
A trajectory may therefore pass Prompt Verification even when its mathematical formulation differs from the reference or the solver terminates under a non-optimal stopping condition.

For a set of evaluated prompts $\mathcal{P}$ and metric
$q\in\{\mathrm{form},\mathrm{solve},\mathrm{opt},\mathrm{verify}\}$,
the corresponding success rate is
\begin{equation}
    R^{q}
    =
    \frac{1}{|\mathcal{P}|}
    \sum_{i\in\mathcal{P}} M_i^{q}.
    \label{eq:metric_rate}
\end{equation}
The same definition is applied to subsets of prompts containing one, two, or three constraint types.
Accordingly, an entry such as $29/30$ in Table~\ref{tab:results} indicates that 29 of the 30 prompts in that subset received a value of one for the corresponding metric.

\subsection{Experimental Results}
Overall, PLUTO performed well across all metrics. Although some runs showed deviations in Constraint Formulation or terminated non-optimally, these did not prevent the tool from producing trajectories that satisfied the prompt in most cases. This indicates that PLUTO and its underlying optimization framework are robust to imperfections at intermediate stages of the pipeline, allowing the tool to succeed even when individual layers are challenged.

\begin{table}[t]
\centering
\caption{Evaluation results by number of constraints per prompt.}
\label{tab:results}
\small
\setlength{\tabcolsep}{4pt}
\renewcommand{\arraystretch}{1.4}
\begin{tabular}{>{\raggedright\arraybackslash}m{2.6cm}cccc}
\hline
 & \multicolumn{3}{c}{Constraints per prompt} & \\
\cline{2-4}
 & 1 & 2 & 3 & Total \\
\hline
Prompts evaluated         & 30    & 10    & 10    & 50 \\
\hline
Constraint Formulation       & 28/30 & 9/10  & 8/10  & 45/50 (90\%) \\
Solver Completion       & 30/30 & 10/10 & 10/10 & 50/50 (100\%) \\
Optimal Convergence       & 28/30 & 6/10  & 8/10  & 42/50 (84\%)\\
\hline
Prompt Verification         & 29/30 & 9/10  & 9/10  & 47/50 (94\%) \\
\hline
\end{tabular}
\end{table}

The distribution of metric results across the sets of one-, two-, and three-constraint prompts is shown in Table~\ref{tab:results}. Solver Completion was achieved for all 50 prompts, with no runs failing or reporting infeasibility, so a trajectory was available for evaluation in every case. Constraint Formulation matched the reference in 90\% of prompts overall, decreasing slightly with complexity from 28 of 30 single-constraint prompts to 9 of 10 two-constraint and 8 of 10 three-constraint prompts. Optimal Convergence was the least consistent metric at 84\% overall. Single-constraint prompts converged optimally in 28 of 30 cases and three-constraint prompts in 8 of 10, while two-constraint prompts dropped to 6 of 10, indicating that non-optimal convergence did not scale directly with the number of constraints. It should be noted that outside of the T10 constraint type, all single constraint scenarios did solve optimally, suggesting that at a core level individual basic constraints do not cause convergence issues. Despite these variations, Prompt Verification remained high across all levels of complexity, with 29 of 30 single-constraint prompts and 9 of 10 prompts in both the two- and three-constraint sets passing, for an overall success rate of 94\%. Notably, the drops in formulation and convergence did not translate into comparable drops in verification, reinforcing that PLUTO can produce trajectories that satisfy the prompt even when intermediate stages of the pipeline are imperfect.

\begin{figure}[t!]
    \centering
    \includegraphics[width=.85\linewidth]{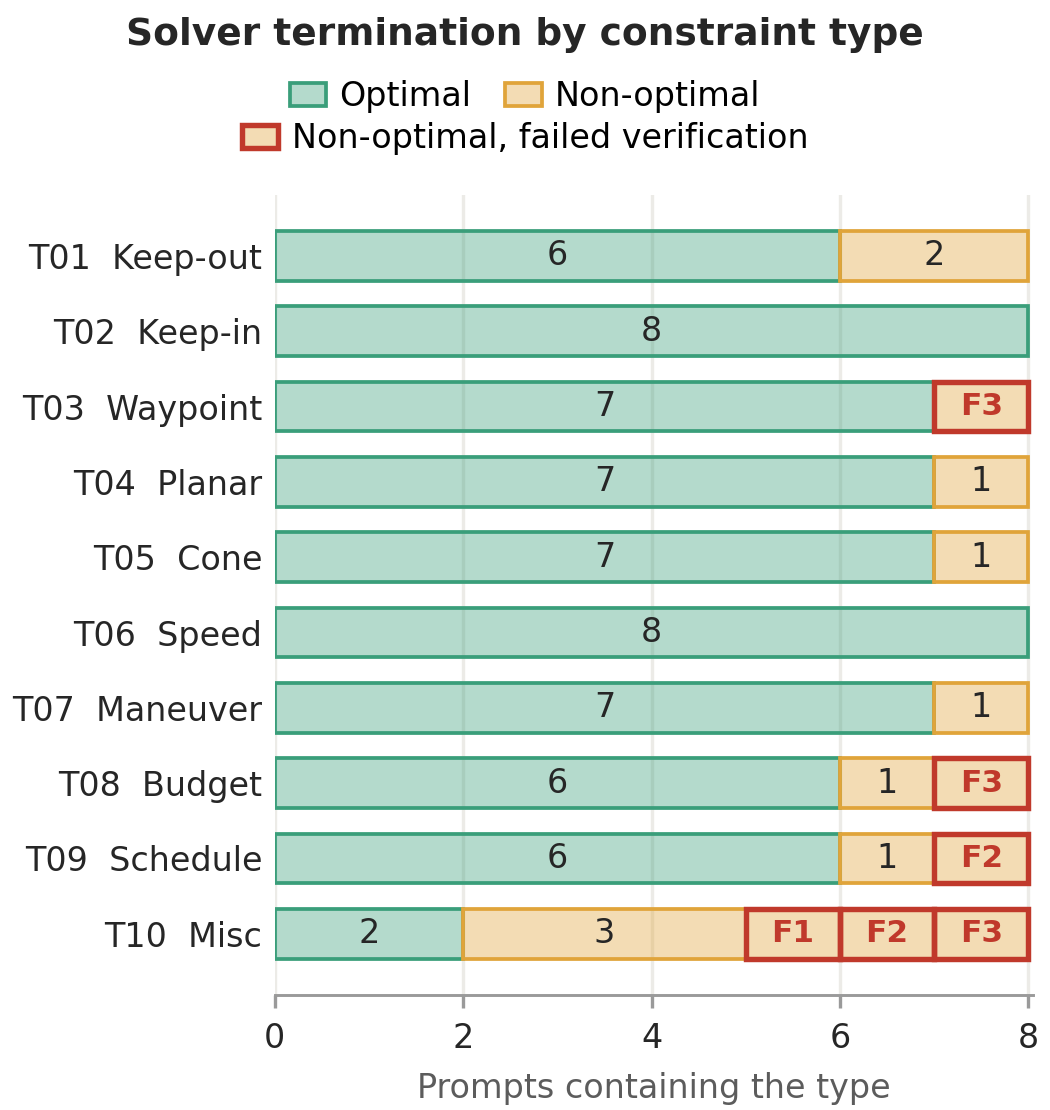}
    \caption{Solver termination by constraint type. Red outlines mark the three runs that failed Prompt Verification (F1--F3).}
    \label{fig:termination}
\end{figure}

Breaking down solver termination by constraint type reveals where convergence issues arise, as shown in Fig.~\ref{fig:termination}. Each bar represents the eight prompts in which that constraint type appears, so a multi-constraint prompt is counted once under each of its types. As a result, the figure shows 80 constraint occurrences rather than the 50 individual prompts. Convex types such as keep-in regions (T02) and speed limits (T06) converged optimally in every prompt, while non-optimal terminations were concentrated in prompts containing T10 constraints, which account for six of the eight non-optimal runs. The non-convex keep-out zone (T01) converged optimally in six of its eight prompts, and its two non-optimal runs still passed verification, indicating that the auto-convexification step handles standard non-convex constraints effectively. Most notably, every run that converged optimally also passed verification, while three of the eight non-optimal runs failed. Optimal Convergence is therefore a strong indicator of success, whereas non-optimal convergence signals that a trajectory warrants closer review rather than indicating failure outright. 

Of the five prompts whose constraint formulations differed from the reference, two still passed Prompt Verification. In the first, a ``corridor'' expected to be modeled as convex planar bounds was instead written as a keep-in cylinder. Because this difference stems from the interpretation of the word ``corridor,'' the agent's output is not incorrect, and the resulting trajectory satisfied the prompt. In the second, the agent's interpretation of the word ``gradual'' changed the implementation but not the verification outcome. Since the agent has no access to the reference or the constraint table, some variation in interpretation is expected, and these cases show that such variation does not always lead to failure. The remaining three deviations, however, did result in verification failures.

\begin{figure*} 
\centering
\captionsetup[subfigure]{font=footnotesize, justification=raggedright,
                         singlelinecheck=false}%
\begin{subfigure}[t]{0.48\textwidth}
    \centering
    \includegraphics[width=\linewidth]{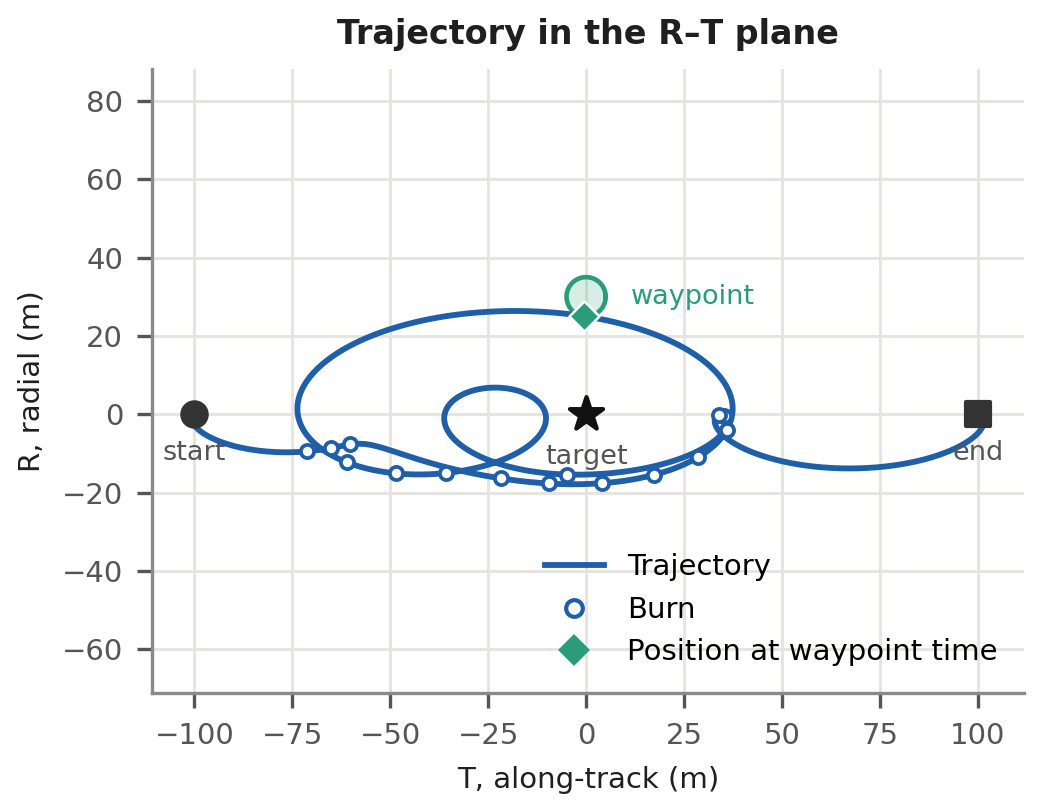}
    \caption{Trajectory plot for the chaser and target.}
    \label{fig:run35_traj}
\end{subfigure}%
\hfill%
\begin{subfigure}[t]{0.48\textwidth}
    \centering
    \includegraphics[width=\linewidth]{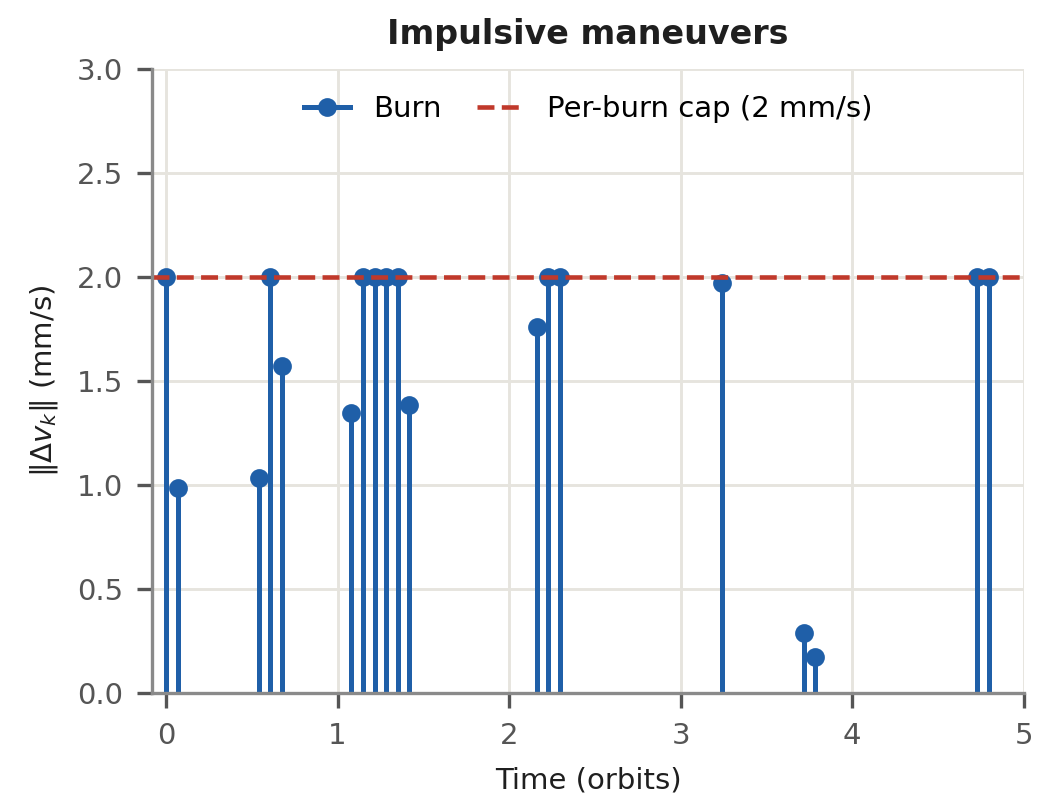}
    \caption{Verification plot for the maneuver limit constraint.}
    \label{fig:run35_verification}
\end{subfigure}
\caption{Outputs for a two-constraint prompt: \textit{``Create a new scenario where the chaser spacecraft passes within 5~m of a point 30~m above the target one third of the way through the trajectory and then passes by the target, with no single maneuver using more than 2~mm/s of delta-v.''}}
\label{fig:example_trajectories}
\end{figure*}

Although many trajectories passed verification despite either a formulation mismatch or non-optimal convergence, all three Prompt Verification failures exhibited both, suggesting that failure arose from compounding issues rather than a single deviation. Another angle to consider is that as shown in Fig.~\ref{fig:termination}, every failure also involved a T10 constraint, indicating that complex or uncommon constraints remain the most challenging for the agent to formulate correctly. The failures were also spread evenly across complexity levels, with one each among the single-, two-, and three-constraint prompts, indicating that they were driven by the constraint type rather than the number of constraints.

Two of the three failures involved the same T10 constraint, an out-of-plane exclusion, and resulted from the same formulation error. The prompt stated that the chaser should remain on \textit{either} side of an exclusion plane. The intended interpretation was to select one side and impose a single convex half-space constraint, a simplification that satisfies the requirement while keeping the problem convex. In both cases, the agent instead attempted to encode the full disjunction, producing a non-convex formulation that led to non-optimal convergence and failed verification. In one of these runs, the agent identified the simpler one-sided formulation but set it aside in favor of a more general solution that preserved both options. The remaining failure followed a similar pattern. The agent began with a correct constraint, but when it encountered a blocker, it added complexity until the formulation became incorrect. Selecting a better free parameter after the initial failure would likely have succeeded.

These failures highlight a core challenge in applying agents to fields that rely heavily on domain expertise. A central part of aerospace engineering practice is recognizing when a simpler formulation is sufficient to meet a requirement, rather than over-engineering a solution to the problem. In these cases, the agent favored completeness over sufficiency, overlooking valid simplifications that an experienced engineer would likely have made.

Much of this judgment comes from experience, which the agent did not have during testing. All prompts were run on a clean version of PLUTO's directory with no stored context from previous runs, ensuring that each prompt was evaluated on equal footing but also leaving the agent without accumulated guidance on simplifications like those described above. A benefit of PLUTO's context-rich structure is that this kind of ``tools of the trade'' knowledge can be encoded into its reference materials, making such shortcuts more accessible to the agent. As tips and lessons from previous runs accumulate with continued use, failures of this kind are expected to become less frequent. Over time, a growing library of formulation strategies could also lead to new and more efficient ways to represent constraints and desired rendezvous behaviors.

Some aspects of the results are harder to capture in numbers, particularly PLUTO's visual outputs. Fig.~\ref{fig:example_trajectories} shows the output plots for a two-constraint prompt combining a waypoint constraint with a per-maneuver delta-v limit. The trajectory in Fig.~\ref{fig:run35_traj} illustrates the nontrivial path required to reach the designated waypoint while keeping each individual maneuver within the delta-v limit, and Fig.~\ref{fig:run35_verification} confirms that the limit is satisfied. These plots keep the user in the loop throughout the design process by showing directly how the intent of each prompt was interpreted and represented in the trajectory. Making these outputs a core part of the tool ensures that, while the agent accelerates constraint formulation, the user retains control over the design and can verify that each trajectory reflects their intent. 

\section{Conclusion}

This paper introduces PLUTO, an agentic framework that integrates an AI coding tool with an SCP backend to enable language-driven trajectory design. Rather than asking the agent to solve the optimization problem directly, PLUTO hones the agent's capabilities to interpret natural-language requirements in a problem formulation space, synthesizing structured constraints that are subsequently classified, convexified when necessary, and solved within a well-defined optimization backend. 
The proposed framework therefore reduces the formulation burden on the user without relinquishing the guarantees provided by numerical optimization and explicit verification.
On the Prompt Verification metric, a 94\% success rate indicates that this separation enables reliable translation of high-level user intent into feasible trajectory design formulations while remaining robust to imperfections in intermediate agent outputs. The remaining challenges primarily concern higher-level engineering judgment, such as selecting the simplest formulation sufficient to satisfy a requirement. Future work will investigate how lessons from prior runs can be incorporated into the framework towards recursive self-improvement (RSI) \cite{xia2026rrsiregularizedrecursiveselfimprovement,lee2026metaharnessendtoendoptimizationmodel} and evaluate its model-agnostic capability using open-source agents.

\section{Acknowledgments}

Large language models were used during this work to assist with writing and debugging code for the PLUTO framework, as well as to refine the language and clarity of the manuscript. All AI-assisted content was carefully reviewed and validated by the authors, who accept full responsibility for the accuracy and integrity of this paper.


\bibliography{main.bib}

@String { proc_IEEE_AC              = {{IEEE Aerospace Conference}} }

@String { proc_IEEE_CDC             = {{Proc.\ IEEE Conf.\ on Decision and Control}} }

@ARTICLE{malyuta_scp_2022,
  author={Malyuta, D. and Reynolds, T. P. and Szmuk, M. and Lew, T. and Bonalli, R. and Pavone, M. and Açıkmeşe, B.},
  journal={IEEE Control Systems Magazine}, 
  title={{Convex Optimization for Trajectory Generation: A Tutorial on Generating Dynamically Feasible Trajectories Reliably and Efficiently}}, 
  year={2022},
  volume={42},
  number={5},
  pages={40-113},
  doi={10.1109/MCS.2022.3187542}}

@INPROCEEDINGS{Banerjee_2020,
  author={Banerjee, Somrita and Lew, Thomas and Bonalli, Riccardo and Alfaadhel, Abdulaziz and Alomar, Ibrahim Abdulaziz and Shageer, Hesham M and Pavone, Marco},
  booktitle=proc_IEEE_AC, 
  title={{Learning-based Warm-Starting for Fast Sequential Convex Programming and Trajectory Optimization}}, 
  year={2020},
  volume={},
  number={},
  pages={1-8},
  doi={10.1109/AERO47225.2020.9172293}}

@phdthesis{damico_phd_2010,
  type = {{PhD} {Thesis}},
  title = {Autonomous {Formation} {Flying} in {Low} {Earth} {Orbit}.},
  language = {English},
  school = {Delft University},
  author = {D'Amico, S.},
  year = {2010}
}

@inproceedings{bonalli_2019_gusto,
    author = {Bonalli, R. and Cauligi, A. and Bylard, A. and Pavone, M.},
    title = {{GuSTO: Guaranteed sequential trajectory optimization via sequential convex programming}},
    booktitle = {IEEE International Conference on Robotics and Automation (ICRA)},
    address = {Montreal, May},
    year = {2019},
}

@article{koenig2017new,
  title={New state transition matrices for spacecraft relative motion in perturbed orbits},
  author={Koenig, Adam W and Guffanti, Tommaso and D’Amico, Simone},
  journal={Journal of Guidance, Control, and Dynamics},
  volume={40},
  number={7},
  pages={1749--1768},
  year={2017},
  publisher={American Institute of Aeronautics and Astronautics}
}

@INPROCEEDINGS{art_ieeeaero24,
  author={Guffanti, Tommaso and Gammelli, Daniele and D’Amico, Simone and Pavone, Marco},
  booktitle=proc_IEEE_AC, 
  title={Transformers for Trajectory Optimization with Application to Spacecraft Rendezvous}, 
  year={2024},
  volume={},
  number={},
  pages={1-13},
  doi={10.1109/AERO58975.2024.10521334}}

@inproceedings{oguri2023successive,
  title={Successive convexification with feasibility guarantee via augmented lagrangian for non-convex optimal control problems},
  author={Oguri, Kenshiro},
  booktitle=proc_IEEE_CDC,
  pages={3296--3302},
  year={2023},
  organization={IEEE}
}

@article{briden2025transformer,
  title={Transformer-Based Tight Constraint Prediction for Efficient Powered Descent Guidance},
  author={Briden, Julia and Gurga, Trey and Johnson, Breanna and Cauligi, Abhishek and Linares, Richard},
  journal={Journal of Guidance, Control, and Dynamics},
  volume={48},
  number={5},
  pages={1054--1070},
  year={2025},
  publisher={American Institute of Aeronautics and Astronautics}
}

@inproceedings{mao2016successive,
  title={Successive convexification of non-convex optimal control problems and its convergence properties},
  author={Mao, Yuanqi and Szmuk, Michael and A{\c{c}}{\i}kme{\c{s}}e, Beh{\c{c}}et},
  booktitle={2016 IEEE 55th Conference on Decision and Control (CDC)},
  pages={3636--3641},
  year={2016},
  organization={IEEE}
}

@inproceedings{kim2022guided,
  title={Guided policy search using sequential convex programming for initialization of trajectory optimization algorithms},
  author={Kim, Taewan and Elango, Purnanand and Malyuta, Danylo and A{\c{c}}{\i}kme{\c{s}}e, Beh{\c{c}}et},
  booktitle={2022 American Control Conference (ACC)},
  pages={3572--3578},
  year={2022},
  organization={IEEE}
}

@inproceedings{foutter2025space-llava,
  author={Foutter, Matthew and Gammelli, Daniele and Kruger, Justin and Foss, Ethan and Bhoj, Praneet and Guffanti, Tommaso and D'Amico, Simone and Pavone, Marco},
  booktitle={2025 IEEE Aerospace Conference}, 
  title={Space-LLaVA: A Vision-Language Model Adapted to Extraterrestrial Applications}, 
  year={2025},
  volume={},
  number={},
  pages={1-23},
  doi={10.1109/AERO63441.2025.11068757}}

@article{takubo2026language,
  title={Language-conditioned safe trajectory generation for spacecraft rendezvous},
  author={Takubo, Yuji and Dwivedi, Arpit and Ramkumar, Sukeerth and Pabon, Luis A and Gammelli, Daniele and Pavone, Marco and D’Amico, Simone},
  journal={Journal of Guidance, Control, and Dynamics},
  pages={1--20},
  year={2026},
  publisher={American Institute of Aeronautics and Astronautics}
}

@article{peng2026nc2c,
  title={NC2C: Automated Convexification of Generic Non-Convex Optimization Problems},
  author={Peng, Xinyue and Liu, Yanming and Cang, Yihan and Zhang, Yuwei and Wang, Xinyi and Deng, Songhang and Cao, Jiannan},
  journal={arXiv preprint arXiv:2601.04789},
  year={2026}
}

@article{jain2026autonomous,
  title={Autonomous Reasoning for Spacecraft Control: A Large Language Model Framework with Group Relative Policy Optimization},
  author={Jain, Amit and Linares, Richard},
  journal={arXiv preprint arXiv:2601.04334},
  year={2026}
}

@inproceedings{ramamonjison2023nl4opt,
  title={Nl4opt competition: Formulating optimization problems based on their natural language descriptions},
  author={Ramamonjison, Rindranirina and Yu, Timothy and Li, Raymond and Li, Haley and Carenini, Giuseppe and Ghaddar, Bissan and He, Shiqi and Mostajabdaveh, Mahdi and Banitalebi-Dehkordi, Amin and Zhou, Zirui and others},
  booktitle={NeurIPS 2022 competition track},
  pages={189--203},
  year={2023},
  organization={PMLR}
}

@article{van2024vernacopter,
  title={VernaCopter: Disambiguated natural-language-driven robot via formal specifications},
  author={van de Laar, Teun and Zhang, Zengjie and Qi, Shuhao and Haesaert, Sofie and Sun, Zhiyong},
  journal={arXiv preprint arXiv:2409.09536},
  year={2024}
}

@inproceedings{wang2024chatstl,
  title={Chatstl: A framework of translation from natural language to signal temporal logic specifications for autonomous vehicle navigation out of blocked scenarios},
  author={Wang, Yujin and Huang, Zhaoyan and Dong, Shiying and Chu, Hongqing and Yin, Xiang and Gao, Bingzhao},
  booktitle={2024 16th International Conference on Computer and Automation Engineering (ICCAE)},
  pages={483--487},
  year={2024},
  organization={IEEE}
}

@inproceedings{chen2023nl2tl,
  title={Nl2tl: Transforming natural languages to temporal logics using large language models},
  author={Chen, Yongchao and Gandhi, Rujul and Zhang, Yang and Fan, Chuchu},
  booktitle={Proceedings of the 2023 Conference on Empirical Methods in Natural Language Processing},
  pages={15880--15903},
  year={2023}
}

@inproceedings{yang2024plug,
  title={Plug in the safety chip: Enforcing constraints for llm-driven robot agents},
  author={Yang, Ziyi and Raman, Shreyas S and Shah, Ankit and Tellex, Stefanie},
  booktitle={2024 IEEE International Conference on Robotics and Automation (ICRA)},
  pages={14435--14442},
  year={2024},
  organization={IEEE}
}

@article{takubo2026orbittamp,
  title={Orbit{TAMP}: Grounding Language Models for Task and Motion Planning in Spacecraft Rendezvous},
  author={Takubo, Yuji and Gammelli, Daniele and Pavone, Marco and D'Amico, Simone},
  journal={arxiv preprint: arXiv:2610.01093v1},
  year={2026}
}

@inproceedings{ahmaditeshnizi2024optimus,
  title     = {{OptiMUS}: Scalable Optimization Modeling with
               ({MI}){LP} Solvers and Large Language Models},
  author    = {Ahmaditeshnizi, Ali and Gao, Wenzhi and Udell, Madeleine},
  booktitle = {Proceedings of the 41st International Conference
               on Machine Learning},
  series    = {Proceedings of Machine Learning Research},
  volume    = {235},
  pages     = {577--596},
  year      = {2024},
  publisher = {PMLR},
  url       = {https://proceedings.mlr.press/v235/ahmaditeshnizi24a.html}
}

@article{feng2025words,
  title={From words to safety: Language-conditioned safety filtering for robot navigation},
  author={Feng, Zeyuan and Zhang, Haimingyue and Bansal, Somil},
  journal={arXiv preprint arXiv:2511.05889},
  year={2025}
}

@inproceedings{ismail2024narrate,
  title={Narrate: Versatile language architecture for optimal control in robotics},
  author={Ismail, Seif and Arbues, Antonio and Cotterell, Ryan and Zurbr{\"u}gg, Ren{\'e} and Alonso, Carmen Amo},
  booktitle={2024 IEEE/RSJ International Conference on Intelligent Robots and Systems (IROS)},
  pages={9628--9635},
  year={2024},
  organization={IEEE}
}

@inproceedings{liang2023code,
  title={Code as policies: Language model programs for embodied control},
  author={Liang, Jacky and Huang, Wenlong and Xia, Fei and Xu, Peng and Hausman, Karol and Ichter, Brian and Florence, Pete and Zeng, Andy},
  booktitle={2023 IEEE International conference on robotics and automation (ICRA)},
  pages={9493--9500},
  year={2023},
  organization={IEEE}
}

@article{zhang2025llm,
title = {A Large Language Model-Based Approach to Spacecraft Formation Mission Planning},
journal = {IFAC-PapersOnLine},
volume = {59},
number = {20},
pages = {2166-2170},
year = {2025},
note = {23rd IFAC Symposium on Automatic Control in Aerospace ACA 2025},
issn = {2405-8963},
doi = {https://doi.org/10.1016/j.ifacol.2025.11.480},
author = {Yuhang Zhang and Bohan Jiao and Zhaohui Dang},
}

@article{rodriguez2024language,
  title={Language models are spacecraft operators},
  author={Rodriguez-Fernandez, Victor and Carrasco, Alejandro and Cheng, Jason and Scharf, Eli and Siew, Peng Mun and Linares, Richard},
  journal={arXiv preprint arXiv:2404.00413},
  year={2024}
}

@misc{xia2026rrsiregularizedrecursiveselfimprovement,
      title={RRSI: Regularized Recursive Self-Improvement of Agent Harnesses}, 
      author={Peng Xia and Rujun Han and Zifeng Wang and Yanfei Chen and Yufan Zhuang and Yoonho Lee and Chengsong Huang and Han Yu and Zhongying CuiZhu and Yifei Ming and Huaxiu Yao and Burak Gokturk and Tomas Pfister and Chen-Yu Lee},
      year={2026},
      eprint={2609.24972},
      archivePrefix={arXiv},
      primaryClass={cs.LG},
      url={https://arxiv.org/abs/2609.24972}, 
}

@misc{lee2026metaharnessendtoendoptimizationmodel,
      title={Meta-Harness: End-to-End Optimization of Model Harnesses}, 
      author={Yoonho Lee and Roshen Nair and Qizheng Zhang and Kangwook Lee and Omar Khattab and Chelsea Finn},
      year={2026},
      eprint={2603.28052},
      archivePrefix={arXiv},
      primaryClass={cs.AI},
      url={https://arxiv.org/abs/2603.28052}, 
}


\thebiography

\begin{biographywithpic}
{Eleanor Brosius}{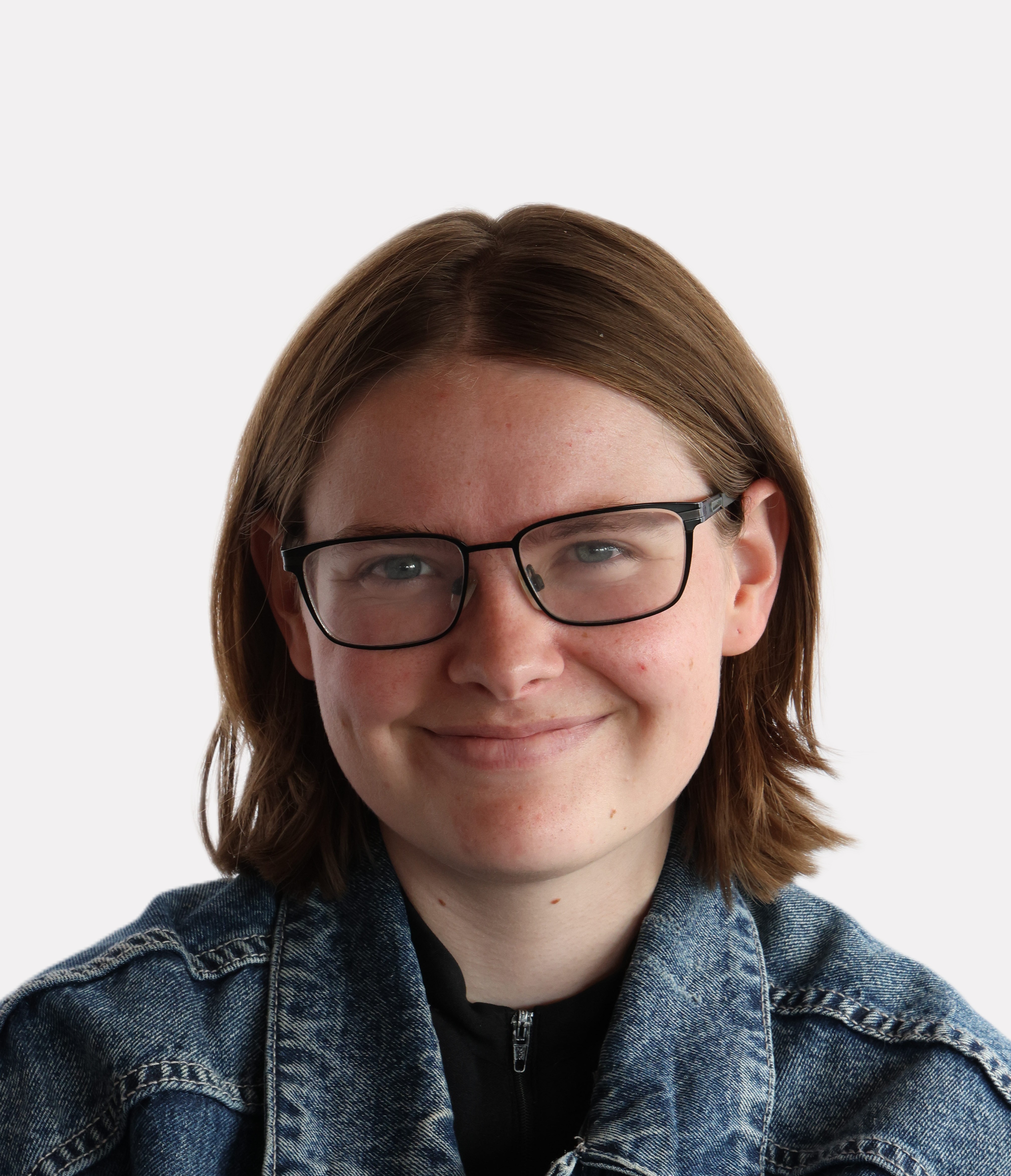} is a master's student in the Department of Aeronautics and Astronautics at Stanford University, where she works in the Autonomous Systems Laboratory. Her research focuses on trajectory optimization for space mission applications and the use of agentic AI to make advanced optimization methods more accessible in aerospace design. She received her B.S. in Electrical Engineering from the University of Washington in 2025, where she worked on CubeSats and autonomous underwater vehicles. 
\end{biographywithpic}%

\begin{biographywithpic}
{Yuji Takubo}{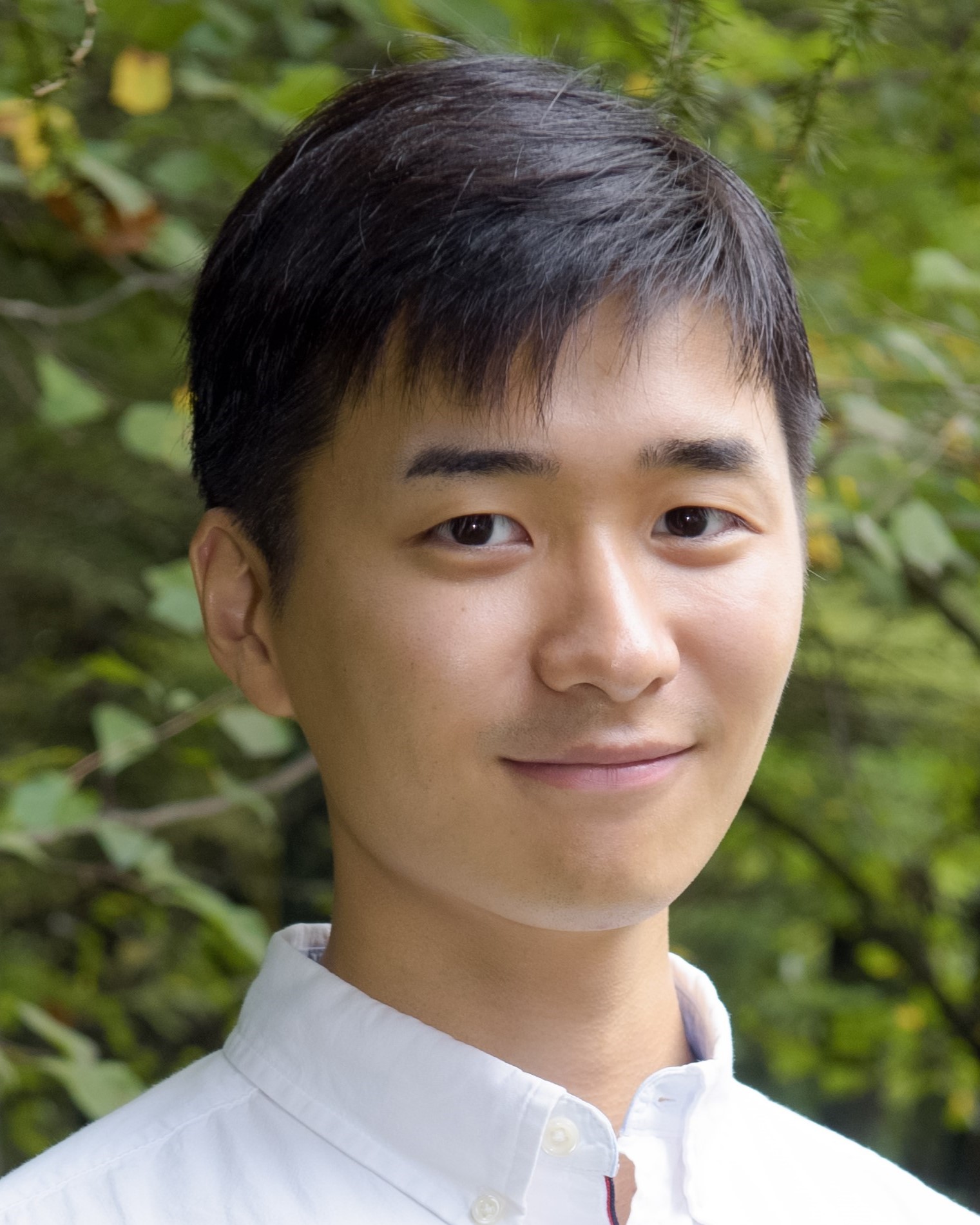} is a Ph.D. candidate in the Space Rendezvous Lab at Stanford University. He received the M.S. in Aeronautics and Astronautics from Stanford University (2025) and the B.S. in Aerospace Engineering from the Georgia Institute of Technology (2023).
He is a recipient of the NASA Jet Propulsion Laboratory Graduate Fellowship (2023) and the Ezoe Memorial Recruit Foundation Fellowship (2020-).
His research focuses on optimization and astrodynamics, particularly with emphasis on distributed space systems.
\end{biographywithpic}%

\begin{biographywithpic}
{Daniele Gammelli}{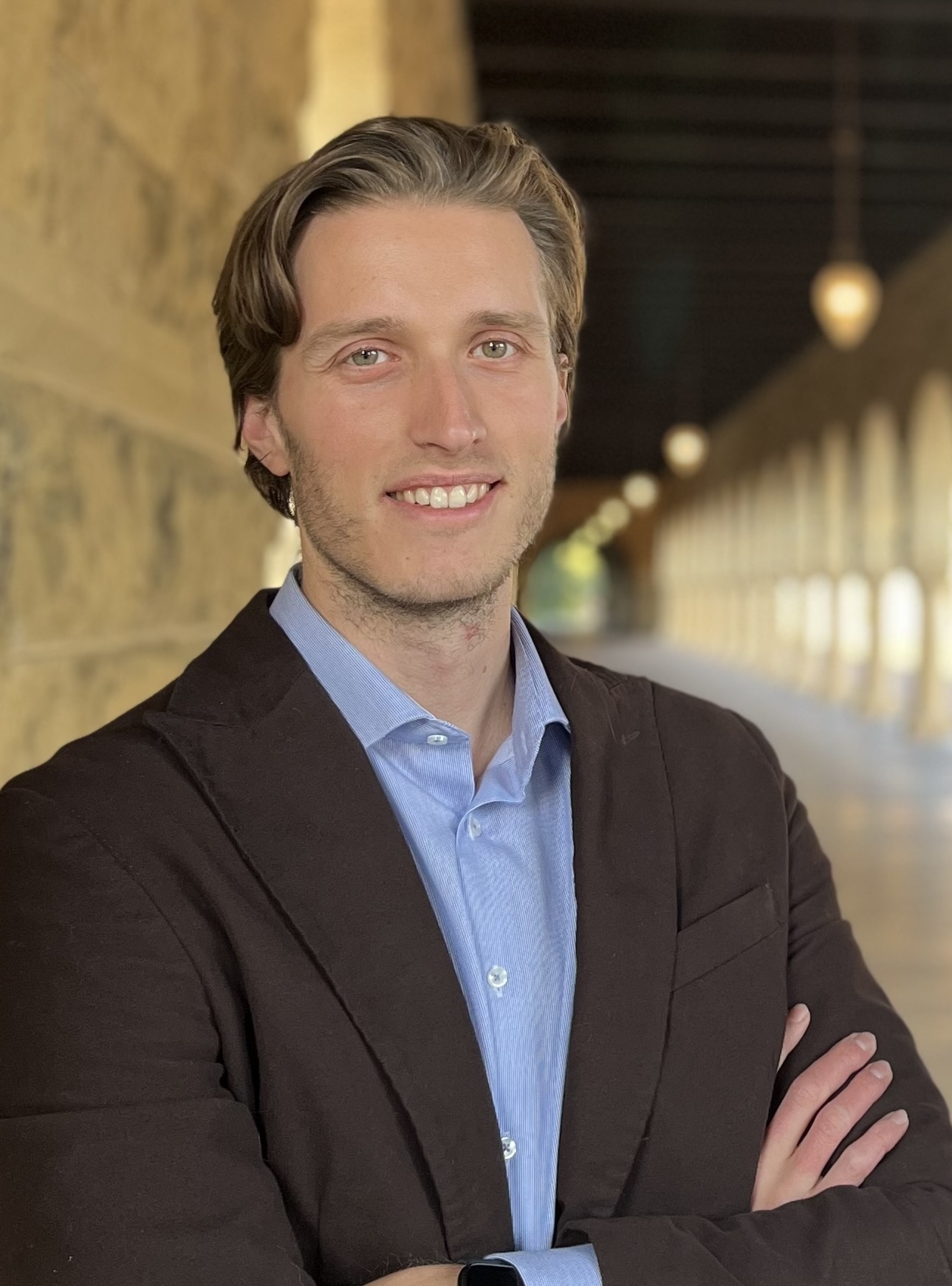} is the Research Director of the Machine Intelligence for Robot Autonomy Laboratory at the Italian Institute of Artificial Intelligence (AI4I), and a Researcher in the Department of Aeronautics and Astronautics at Stanford University. Since 2022, he has been a Research Fellow at the Center for Aerospace Autonomy Research (CAESAR) at Stanford and, until 2025, at the Center for Automotive Research at Stanford (CARS). He received his Ph.D. in Machine Learning and Mathematical Optimization from the Technical University of Denmark (DTU) in 2022, where his doctoral thesis was nominated for the DTU Best Thesis of the Year Award and the DTU Young Researcher Award. Dr. Gammelli’s research focuses on developing the algorithmic foundations and system-level methodologies that enable AI-powered autonomous systems to operate safely, efficiently, and reliably in high-stakes environments, with an emphasis on aerospace autonomy and next-generation mobility systems.
\end{biographywithpic}

\begin{biographywithpic}
{Simone D'Amico}{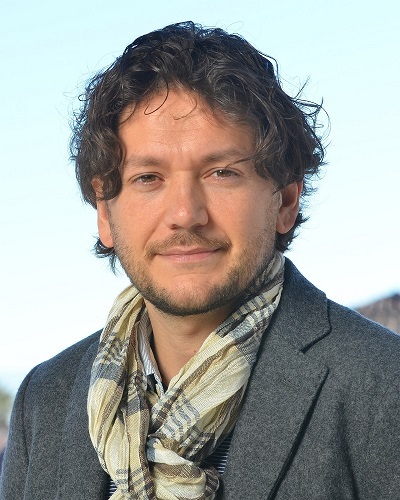} (Member, IEEE) is an Associate Professor of Aeronautics and Astronautics (AA) and Professor of Geophysics (by Courtesy). He is the Founding Director of the Stanford Space Rendezvous Laboratory, Director of the Undergraduate Program in AA, and Founding Co-Director of the Center for AEroSpace Autonomy Research (CAESAR) at Stanford. He is a W.M. Keck Faculty Scholar in the School of Engineering and Hoover Science Fellow.
He has 20+ years of experience in research and development of autonomous spacecraft and distributed space systems, including multi-agent architectures such as “rendezvous, proximity operations, and capture”, “satellite formation-flying and swarms”, “fractionated spacecraft”, and “mega-constellations”. He developed the GNC system of several such missions (e.g., TanDEM-X, PRISMA) and is currently the institutional PI of four autonomous satellite swarms funded by NASA and NSF with one of them (NASA Starling) operational in orbit.
Besides academia, Dr. D’Amico is the Chief Science Officer (CSO) of EraDrive. He is on the Advisory Board of four space start-ups focusing on distributed space systems for future applications in SAR remote sensing, orbital lifetime prolongation, and space-based solar power. He received the B.S. and M.S. degrees from Politecnico di Milano (2003) and the Ph.D. degree from Delft University of Technology (2010). Before Stanford, Dr. D’Amico was research scientist and team leader at the German Aerospace Center (DLR). He was the recipient of several awards, including the 2020 IEEE M. Barry Carlton Award, Best Paper Awards at AMOS (2025), IAF (2022), IEEE (2021), AIAA (2021), AAS (2019) conferences, the Leonardo 500 Award by the Leonardo da Vinci Society/ISSNAF (2019), FAI/NAA’s Group Diploma of Honor (2018), DLR’s Sabbatical/Forschungssemester (2012) and Wissenschaft Preis (2006), and NASA’s Group Achievement Awards for the STARLING (2024) and GRACE mission (2004).
\end{biographywithpic}

\begin{biographywithpic}
{Marco Pavone}{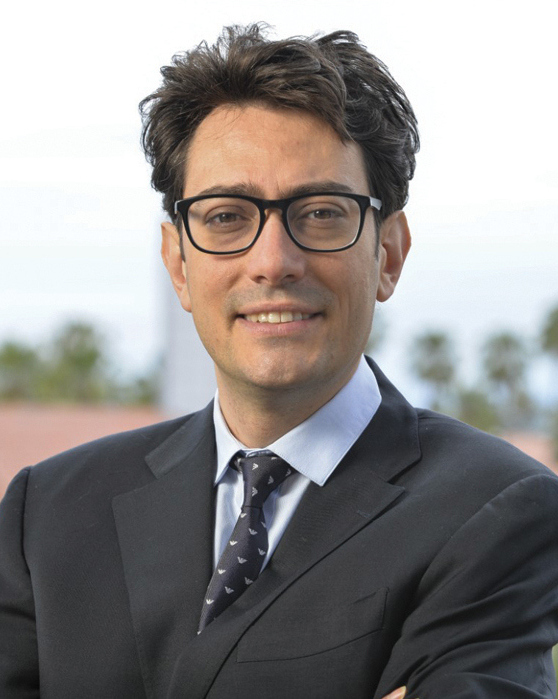} (Member, IEEE) is a Professor of Aeronautics and Astronautics at Stanford University, where he directs the Autonomous Systems Laboratory and the Center for Automotive Research at Stanford. He is also a Distinguished Research Scientist at NVIDIA where he leads autonomous vehicle research. Before joining Stanford, he was a Research Technologist within the Robotics Section at the NASA Jet Propulsion Laboratory. He received a Ph.D. degree in Aeronautics and Astronautics from the Massachusetts Institute of Technology in 2010. His main research interests are in the development of methodologies for the analysis, design, and control of autonomous systems, with an emphasis on self-driving cars, autonomous aerospace vehicles, and future mobility systems. He is a recipient of a number of awards, including a Presidential Early Career Award for Scientists and Engineers from President Barack Obama, an Office of Naval Research Young Investigator Award, a National Science Foundation Early Career (CAREER) Award, a NASA Early Career Faculty Award, and an Early-Career Spotlight Award from the Robotics Science and Systems Foundation. He was identified by the American Society for Engineering Education (ASEE) as one of America's 20 most highly promising investigators under the age of 40. His work has been recognized with best paper nominations or awards at a number of venues, including the European Conference on Computer Vision, the IEEE International Conference on Robotics and Automation, the European Control Conference, the IEEE International Conference on Intelligent Transportation Systems, the Field and Service Robotics Conference, the Robotics: Science and Systems Conference, and the INFORMS Annual Meeting.
\end{biographywithpic}

\end{document}